\documentclass[12pt]{article}
\usepackage[cp1251]{inputenc}
\usepackage[russian]{babel}
\begin{document}

\author{Beloshapka V.K.}

\date{}

\title{Associative algebras in $CR$-geometry}

\maketitle

\begin{abstract}
A procedure for the algebraization of a $CR$-manifold and its holomorphic automorphisms is described. Examples of the application of algebraization are considered. Questions arising in connection with the algebraization of a $CR$-manifold are formulated. The possibilities of extending this procedure to other branches of geometry and analysis are discussed.
 \end{abstract}

\footnote{
Mechanics and Mathematics Faculty, Lomonosov Moscow State University, Vorobyovy Gory, 119992 Moscow, Russia, vkb@strogino.ru }

{\bf 1. Introduction}

 \vspace{3ex}

$CR$-geometry studies the geometric and analytic properties of real submanifolds of complex manifolds up to holomorphic equivalence. In $CR$-geometry, the method of a model surface has been formed. Its essence is that the germ of a real submanifold is associated with a simpler object -- its model surface, which is real algebraic and whose properties strongly reflect the properties of the original germ. In particular, the local group of holomorphic automorphisms of the germ is majorized in dimension by the group of the model surface. The set of model surfaces of a fixed type is a finite-dimensional space. For a generic point of this space, i.e., the corresponding model surface, the group of holomorphic automorphisms is predictably small. Model surfaces with a large nontrivial group exist, but are infrequent and of great interest. In this context, procedures that allow multiplying the number of such examples are of value. The simplest way is to construct the Cartesian product of model surfaces.

In this paper, a much more general procedure will be considered -- $\mathcal{S}$-algebraization. This procedure is based on replacing the field of real numbers $\mathbf{R}$ with an arbitrary real finite-dimensional associative commutative algebra with identity $\mathcal{S}$. In this case, accordingly, instead of the field of complex numbers $\mathbf{C}$, we use $\mathcal{S}_c$ -- the complexification of the algebra $\mathcal{S}$. Note that in algebra, algebraic geometry, and algebraic number theory, procedures of this kind (extensions and restrictions of the numerical universe) have been used effectively for a long time.

In this paper, by an {\it algebra} we will mean exactly such an algebra $\mathcal{S}$: real finite-dimensional associative commutative with identity (not to be confused with Lie algebras, which are also present in this paper). The description of all such algebras is given by the Wedderburn-Artin theorem (1908, 1927, see \cite{A50},\cite{AT},\cite{LT}). According to this theorem, $\mathcal{S}$ is a direct sum of a finite set of local algebras $\mathcal{S}_j$ (spectrum), i.e., algebras with a unique maximal ideal $m_j$. The quotient $\mathcal{S}_j / m_j$ is a field. This is either $\mathbf{R}$ (a real point of the spectrum) or $\mathbf{C}$ (a complex point of the spectrum). The ideal $m_j$ is nilpotent, it consists of all nilpotent and therefore non-invertible elements, and has codimension one or two in $\mathcal{S}_j$. The set of algebras of dimension no higher than three (up to isomorphism) is finite. In dimension four, parameters appear. The description of all algebras of dimension six and higher seems to be a very difficult problem.

It is not difficult to show that each such algebra can be realized as a subalgebra of a matrix algebra or as a quotient algebra of a polynomial algebra in several variables. Based on the realization as a matrix subalgebra, we can endow $\mathcal{S}$ with a norm that introduces a topology. At the same time, we can require this norm to be {\it submultiplicative}, i.e., $\|A\,B\|\leq \|A\|\,\|B\|$. All operator norms satisfy this condition. Here is their description. We will consider the matrix algebra as a coordinate notation for the algebra of linear operators on some real linear space. Let us fix an arbitrary norm on this vector space. This choice induces a norm on the algebra of linear operators. Namely, the maximum of the norms of the operator on the unit ball. Fixing this norm as the matrix norm, we obtain a submultiplicative norm on $\mathcal{S}$. There are other ways to define a submultiplicative norm on $\mathcal{S}$.

Now we can consider convergent power series on $\mathcal{S}$, both in one and in several variables. The technique is the same as with power series in $\mathbf{R}$ or $\mathbf{C}$. Geometric progression, Abel's lemma, domains and polydisks of convergence. As a special feature, the following should be noted. Any element of the complexified algebra $\mathcal{S}_c$ according to the Wedderburn-Artin theorem is a sum of the form $z+\nu$, where $z$ is an element of $\mathbf{C}^N$, and $\nu$ is a nilpotent, i.e., ${\nu}^m=0$ for some $m$. Writing the Taylor expansion of our series at the point $z$ in powers of $\nu$, we see that the series terminates. Thus, only the scalar coordinates $z$ participate in the description of the domain of convergence, while the nilpotent ones do not.

 \vspace{3ex}

{\bf Statement 1:} Let
$\sum_{\alpha_j \geq 0} \, C_{\alpha_1, \dots ,\alpha_N} \, Z_1^{\alpha_1} \dots Z_N^{\alpha_N}$, where $C_{\alpha_1, \dots,\alpha_N}, Z_j^{\alpha_j} \in \mathcal{S}_c$ be a power series converging in a neighborhood of the origin in $\mathcal{S}^N_c$. Its sum is identically equal to zero if and only if $C_{\alpha_1, \dots ,\alpha_N}=0$ for all $\alpha$.\\
{\it Proof:} Let the series be equal to zero. Consider the subspace of scalars in $\mathcal{S}_c^N$ -- $\mathbf{C}^N$ (all nilpotent coordinates are equal to zero). We obtain a power series in complex variables (with algebraic coefficients)
$\sum \, C_{\alpha_1, \dots ,\alpha_N} \, z_1^{\alpha_1} \dots z_N^{\alpha_N}$. Substituting $z=0$, we obtain $C_{0 \dots  0}=0$. Such a series can be differentiated term by term with respect to any complex variable. Differentiating the required number of times and substituting $z=0$, we obtain $C_{\alpha_1, \dots ,\alpha_N}=0$. The converse is obvious. The statement is proved.

 \vspace{3ex}

From the general description of algebras (finite-dimensional, associative, commutative, with identity), i.e., from the Wedderburn-Artin theorem, there is a simple corollary.

 \vspace{1ex}

{\bf Statement 2:} The set of non-invertible elements of an algebra is a subspace of positive codimension.

\vspace{1ex}

In particular, the identity has a neighborhood consisting of invertible elements.

 \vspace{3ex}

{\bf 2. Algebraization of a real algebraic set}

 \vspace{3ex}

Fix some algebra $\mathcal{S}$. Let $M$ be a nonempty real algebraic subset of $\mathbf{C}^N$. That is, $M$ is the set of common zeros of a finite set of real-valued polynomials on $\mathbf{C}^N\simeq \mathbf{R}^{2N}$.
  \begin{equation}\label{M}
 \{z =x+i\,y \in \mathbf{C}^N: \tilde{\rho}(x,y)=\rho(z,\bar{z})=(\rho_1(z,\bar{z}),\dots,\rho_s(z,\bar{z}))=0\}.
 \end{equation}
Let $\mathcal{S}_c \simeq\mathcal{S}^2$ be the complexification of $\mathcal{S}$. Let further
 $ \bar{Z} =\overline{(X+i \, Y)} =  X-i \, Y, \; X,Y \in \mathcal{S}$ be the complex conjugation defined on $\mathcal{S}_{c}$, and $\mathrm{Re}\, Z=X, \; \mathrm{Im}\, Z=Y$ be the real and imaginary parts.

 \vspace{1ex}

{\it Definition 1:} The $\mathcal{S}$-algebraization of the algebraic set $M$ given by relations (\ref{M}) is the real algebraic set $\mathcal{M}_{\mathcal{S}}=\mathcal{M}$, given by the relations
\begin{equation}\label{SM}
\{Z \in \mathcal{S}_c^N \simeq \mathcal{S}^{2N}: \tilde{\rho}(X,Y)=\rho(Z,\bar{Z})=(\rho_1(Z,\bar{Z}),\dots,\rho_s(Z,\bar{Z}))=0\},
\end{equation}
i.e., it is the result of substituting $(x \rightarrow X, \, y \rightarrow Y, \, z \rightarrow Z, \, \bar{z} \rightarrow \bar{Z})$ into the polynomials defining $M$, i.e., into $\rho_j$.

 \vspace{3ex}

The fact that there is an identity in the algebra $\mathcal{S}$ allows us to identify $(z_1, \dots, z_N) \in \mathbf{C}^N$ and $(z_1 \cdot 1 , \dots, z_N \cdot 1) \in \mathcal{S}_c^N$, to assume that $\mathbf{C}^N \subset \mathcal{S}_c^N$, and that $M \subset \mathcal{M}$ (i.e., $\mathcal{M}$ is not empty).

 \vspace{3ex}

The following statement can be verified directly.

\vspace{1ex}

{\bf Statement 3:} The result of applying two algebraizations to the set $M$; the first with algebra $\mathcal{S}_1$ and the second with algebra $\mathcal{S}_2$ -- is an algebraization with algebra $\mathcal{S}=\mathcal{S}_2 \otimes_{\mathbf{R}} \mathcal{S}_1$.

 \vspace{3ex}

Let us subject an arbitrary nondegenerate model surface $Q$ to algebraization. $"$Arbitrary$"$ means of arbitrary Bloom-Graham type and, thus, of arbitrary $CR$-type $(n,k)$ ($n$ is the dimension of the complex tangent, $k$ is the codimension) in the space $\mathbf{C}^{(n+k)}$ \cite{VB20}. Moreover, it can be a weighted model surface \cite{VB23}. The equations of such a surface have the form
\begin{equation}\label{MQ}
\mathrm{Im}\,w=\Phi(z,\bar{z},\mathrm{Re}\,w),    \;\; \mbox{where  }   w \in \mathbf{C}^k, \;   z \in \mathbf{C}^n,
\end{equation}
where the coordinates of the mapping $\Phi$ are homogeneous or quasi-homogeneous real polynomials with a fixed choice of weights according to the Bloom-Graham-Stepanova type \cite{ST}. Associated with each such surface is a class of germs of real submanifolds in $\mathbf{C}^{(n+k)}$ of the same type $M_0$ with a local equation of the form
\begin{equation}\label{MM}
\mathrm{Im}\,w=\Phi(z,\bar{z},\mathrm{Re}\,w) + \dots,
\end{equation}
where the ellipsis means that in each coordinate the homogeneous polynomial can be perturbed by terms of higher weight.
Holomorphically invariantly associated with each germ are two objects $\mathrm{aut}\,M_0$ and $\mathrm{Aut}\,M_0$ -- the real Lie algebra of holomorphic fields in a neighborhood of zero tangent to $M_0$ and the local group of holomorphic automorphisms of the germ generated by this Lie algebra.

$"$Nondegeneracy$"$ means the fulfillment of two generally independent conditions on $Q$:
  \begin{eqnarray}\label{NDC}
 \nonumber  (*) \quad \mbox{the coordinates of the form }  \; \Phi \;  \mbox{ are linearly independent}, \qquad\qquad\qquad\qquad\\
\nonumber  \mbox{(this guarantees the finiteness of the Bloom-Graham type)},\qquad\qquad\qquad\\
(**) \quad \mathcal{Q} \quad \mbox{is holomorphically nondegenerate (see \cite{BER}).}\qquad\qquad
 \end{eqnarray}
This is the general definition of the nondegeneracy of a model surface; it represents a criterion for the finite dimensionality of ${\rm aut} \, Q$ \cite{BER},\cite{VB20}. At the same time, the dimension of the automorphism group of surface (\ref{MQ}) majorizes the dimension of the automorphism group of germ (\ref{MM}).

 \vspace{3ex}

Thus, the $\mathcal{S}$-algebraization of the model surface $Q$ given by relation (\ref{MQ}) is the model surface $\mathcal{Q}_{\mathcal{S}}=\mathcal{Q}$, given by the relation
\begin{equation}\label{SQ}
\mathrm{Im}\,W=\Phi(Z,\bar{Z},\mathrm{Re}\,W),    \;\; \mbox{where  }   W \in \mathcal{S}_{c}^k, \;   Z \in \mathcal{S}_{c}^n.
\end{equation}

\vspace{3ex}

As already mentioned, the context of $CR$-geometry is the study of real submanifolds of a complex space up to invertible holomorphic transformations. The holomorphy of a function or mapping is equivalent to the possibility of a local representation as a sum of a power series in local complex coordinates. Let there be a holomorphic mapping of the algebra $\mathcal{S}$ into itself in a neighborhood of zero, i.e., $f: \mathcal{S}_c \rightarrow \mathcal{S}_c$ (denote the argument by $Z \in \mathcal{S}_c$). If the dimension of $\mathcal{S}$ is $l$, then this mapping can be viewed as a mapping of $\mathbf{C}^l$ into itself. To identify $\mathcal{S}_c$ and $\mathbf{C}^l$, one needs to choose a basis in $\mathcal{S}$, denote it by $(1,e_2, \dots, e_l)$, then
$$Z=z_1 \, 1+z_2 \,e_2 + \dots +z_l \, e_l=(z_1, \dots,z_l)=z.$$
This identification allows introducing a norm on $\mathcal{S}_c$.
The holomorphy of $f$ in a neighborhood of zero means that $f(Z)=f(z)=\sum_0^{\infty}\, f_j(z)$, where $f_j(z)$ is a homogeneous mapping of $\mathbf{C}^l$ into itself of degree $j$, and the series converges in a neighborhood of zero.

\vspace{1ex}

{\it Definition 2:} We say that $f$ is $\mathcal{S}$-holomorphic in a neighborhood of zero if $f$ is the sum of an $\mathcal{S}$-power series, i.e., $f(Z)= \sum_0^{\infty}\, c_j \, Z^j, \;\; c_j \in \mathcal{S}_c$, converging (in norm) in a neighborhood of zero.

We understand the $\mathcal{S}$-holomorphy of a mapping from $\mathcal{S}^{\alpha}_c$ to $\mathcal{S}^{\beta}_c$ similarly. Thus, $\mathcal{S}$-holomorphy is equivalent to the representability of each coordinate mapping as a convergent power series in several algebraic variables with algebraic coefficients. Or, equivalently, representability as the sum of a convergent series of homogeneous forms in algebraic variables.

\vspace{1ex}

If in the space $\mathcal{S}^N_c \simeq \mathbf{C}^{lN}$ there are two germs -- $M_p$ and $M'_{p'}$, then besides the relation of their holomorphic equivalence, it is appropriate to introduce the concept of $\mathcal{S}$-holomorphic equivalence.

\vspace{1ex}

{\it Definition 3:} We say that the germs of $CR$-manifolds $M_p$ and $M'_{p'}$ are $\mathcal{S}$-{\it holomorphically equivalent} if there exists an $\mathcal{S}$-holomorphic mapping of $M_p$ onto $M'_{p'}$. Notation:
$M_p \approx_{\mathcal{S}} M'_{p'}$.

\vspace{1ex}

Any holomorphic real vector field in a neighborhood of zero in $\mathcal{S}_c^{(n+k)}$ has the form
\begin{equation}\label{VF}
 X=2 \,\mathrm{Re}\, \left(f(Z, W) \,\frac{\partial}{\partial Z} +g(Z,W)\, \frac{\partial}{\partial W}\right),
\end{equation}
where $f$ and $g$ are holomorphic in a neighborhood of zero.

\vspace{1ex}

{\it Definition 4:} We will call a vector field of the form (\ref{VF}) $\mathcal{S}$-holomorphic in a neighborhood of zero, if $f$ and $g$ are $\mathcal{S}$-holomorphic in a neighborhood of zero.

\vspace{3ex}

It is clear that $\mathcal{S}$-holomorphic mappings are a subclass of the class of holomorphic mappings, and $\mathcal{S}$-holomorphic vector fields are a subspace of the space of holomorphic vector fields. Therefore, whenever there is an identification of the ambient space $\mathbf{C}^{l \,N}$ containing our model surface $Q$ with the space $\mathcal{S}_c^{N}$ ($l=\mathrm{dim}\,\mathcal{S}$), we can extract $\mathcal{S}$-holomorphic vector fields in the Lie algebra $\mathrm{aut}\,Q$, and holomorphic automorphisms in $\mathrm{Aut}\,Q$.
Let us denote them by $\mathrm{aut}_{\mathcal{S}}\,Q$ and $\mathrm{Aut}_{\mathcal{S}}\,Q$, respectively.

\vspace{1ex}

Let us have a pair: an algebra $\mathcal{L}$ and its subalgebra $\mathcal{S}$.
It is clear that all conditions in the sense of $\mathcal{S}$ follow from the corresponding conditions in the sense of $\mathcal{L}$.
If we take the minimal possible subalgebra $\mathbf{R}$ as $\mathcal{S}$, then all definitions transform into the usual ones: holomorphy, holomorphic equivalence, etc.

\vspace{3ex}

{\bf Theorem 4:} Let $Q$ be a nondegenerate model surface in $\mathbf{C}^{l\,N} \simeq \mathcal{S}_c^l$, then\\
(a) $\mathrm{aut}_{\mathcal{S}}\,Q$ is a Lie subalgebra in $\mathrm{aut}\,Q$,\\
(b) $\mathrm{Aut}_{\mathcal{S}}\,Q$ is a subgroup in the local group $\mathrm{Aut}\,Q$.\\
{\it Proof:} To prove (a), it is sufficient to prove that the commutator of $\mathcal{S}$-holomorphic vector fields is an $\mathcal{S}$-holomorphic vector field. But this is a consequence of the fact that the result of applying an $\mathcal{S}$-holomorphic $(1,0)$-field to an $\mathcal{S}$-holomorphic function is an $\mathcal{S}$-holomorphic function. Item (a) is proved. Further, $\mathrm{Aut}_{\mathcal{S}}\,Q$ is the image of $\mathrm{aut}_{\mathcal{S}}\,Q$ under the exponential mapping. Thus, to prove (b), it is sufficient to prove that if the right-hand side of a system of first-order ordinary differential equations, resolved with respect to the derivative, is $\mathcal{S}$-holomorphic, then the solution will also be $\mathcal{S}$-holomorphic. This is proved by a standard construction. First, the form of a formal solution is written out as an $\mathcal{S}$-power series, and then (by the method of majorant series) its convergence. The theorem is proved.

\vspace{3ex}

{\bf 3. $\mathrm{RAQ}$-quadrics}

 \vspace{3ex}

Consider the surface
\begin{equation}\label{SP}
Q=\{(z,w) \in \mathbf{C}^2:  \mathrm{Im}\,w=z\,\bar{z}\}.
\end{equation}
$\mathcal{Q}$ is the usual 3-dimensional sphere in $\mathbf{C}^2$ (in a projective realization). The model role of this hypersurface was demonstrated in the first paper on $CR$-geometry \cite{P}. If we assign weights to the variables as follows $[z]=1, \;[w]=2$, then
$\mathrm{aut}\,Q$ consists of five components $g_{-2}+g_{-1}+g_{0}+g_{1}+g_{2}$, with
\begin{eqnarray}\label{ASP}
g_{-2}=2 \, {\rm Re}\,( \alpha \,\frac{\partial}{\partial w}), \; \alpha \in \mathbf{R},\\
\nonumber g_{-1}=2 \, {\rm Re}\,(\beta\, \frac{\partial}{\partial z} +2 \,i \ \bar{\beta}\, z \, \frac{\partial}{\partial w}), \; \beta \in \mathbf{C},\\
\nonumber g_0=2 \, {\rm Re}\,(\gamma \, z \,\frac{\partial}{\partial z} +2 \, {\rm Re} ( \gamma )\,w \, \frac{\partial}{\partial w}), \;  \gamma \in \mathbf{C},\\
\nonumber g_1=2 \, {\rm Re}\,((\delta\, w +2 \,i \, \bar{\delta} \,z^2 )\, \frac{\partial}{\partial z} + 2 \,i \, \bar{\delta} \,z \,w \, \frac{\partial}{\partial w}), \; \delta \in \mathbf{C}.\\
\nonumber g_2=2 \, {\rm Re}\,(\varepsilon \, z\, w \,\frac{\partial}{\partial z} +\varepsilon \,w^2 \, \frac{\partial}{\partial w}), \; \varepsilon \in \mathbf{R}.
\end{eqnarray}
These components generate the transformations
\begin{eqnarray}\label{AUT}
\nonumber \left( z \rightarrow z+p, \quad   w  \rightarrow w +2\,i\, \bar{p}\,z +(q+i\,p\,\bar{p})\right),\\
\left( z \rightarrow \lambda \, z, \quad   w  \rightarrow \lambda \,\bar{\lambda} \,w \right),  \qquad\qquad\qquad\qquad\\
\nonumber  z \rightarrow   \left(1-2\,i\,\bar{a}\,z-(r+i\,a\,\bar{a})   \right)^{-1} \,(z+a w), \\
\nonumber w  \rightarrow   \left(1-2\,i\,\bar{a}\,z-(r+i\,a\,\bar{a})   \right)^{-1} \, w ,
\end{eqnarray}
where $q, \; r \in \mathbf{R}, \;\; p, \; a \in \mathbf{C}, \; \lambda \in \mathbf{C}^{*}$. These transformations, in turn, generate ${\rm Aut} \, Q$, and the dimension of the group is equal to 8.

\vspace{3ex}

In the papers of Ezhov and Schmalz \cite{ES94}, \cite{ES98}, the algebraization of the surface $Q$ was considered, i.e., a surface of the form
\begin{equation}\label{SSP}
 \mathcal{Q}=\{(Z,W) \in \mathcal{S}_{c} \oplus \mathcal{S}_{c}:  {\rm Im}\,W=Z\,\bar{Z} \}.
\end{equation}
$\mathcal{Q}$ is a model surface in $\mathbf{C}^{2l}$, it has $CR$-type $(l,l)$.
Surfaces of the form (\ref{SP}), as well as (\ref{SSP}), are special cases of a general quadratic model surface of arbitrary $CR$-type $(n,k)$ (if the surface is nondegenerate, then $k \leq n^2$). That is, surfaces in $\mathbf{C}^{n+k}$ of the form
\begin{equation}\label{QQ}
 \{  (z,w):     {\rm Im}\,w=<z,\bar{z}>, \; z \in \mathbf{C}^n, \; w \in \mathbf{C}^k\},
\end{equation}
where $<z,\bar{z}>$ is an $\mathbf{R}^k$-valued Hermitian form (see \cite{VB88}).
Surfaces of the form (\ref{SSP}) were called $RAQ$-quadrics (real associative quadrics) in \cite{ES98}, in contrast to surfaces (\ref{QQ}) of a more general form ($CR$-quadrics).

Condition (\ref{NDC}) of nondegeneracy for a $CR$-quadric reduces to the fulfillment of the following conditions
\begin{eqnarray}\label{NDG}
\nonumber  \mbox{ The scalar coordinates of the form  } <z,\bar{z}>        \mbox{ are linearly independent },\\
\mbox{if  }  b \in\mathbf{C}^n     \mbox{ and  }  <z,\bar{b}>=0,        \mbox{ then  } b=0.  \qquad\qquad\qquad\qquad
\end{eqnarray}

\vspace{3ex}

Let us list some results from \cite{ES98}. \\
-- (a) Any $RAQ$-quadric $\mathcal{Q}_{\mathcal{S}}$ is nondegenerate,\\
-- (b) $\mathcal{Q}_{\mathcal{S}_1}$ is holomorphically equivalent to $\mathcal{Q}_{\mathcal{S}_2}$ if and only if $\mathcal{S}_1$ is isomorphic to $\mathcal{S}_2$,\\
-- (c) $\mathcal{Q}_{\mathcal{S}}=\mathcal{Q}_{\mathcal{S}_1} \times \mathcal{Q}_{\mathcal{S}_2}$ (Cartesian product) if and only if $\mathcal{S}=\mathcal{S}_1 \oplus \mathcal{S}_2$ (direct sum).\\
-- (d) The subgroup of the group of holomorphic automorphisms of an $RAQ$-quadric generated by fields of positive weight has the form
\begin{eqnarray}\label{AUT}
 Z \rightarrow   \left(1-2\,i\,\bar{A}\,Z-(R+i\,A\,\bar{A})   \right)^{-1} \,(Z+A W), \\
\nonumber W  \rightarrow   \left(1-2\,i\,\bar{A}\,Z-(R+i\,A\,\bar{A})   \right)^{-1} \, W
\end{eqnarray}
where $R \in \mathcal{S}, \; A \in \mathcal{S}_c$.

\vspace{3ex}
Let us give a comment on item (d). The fact that mappings of the form (\ref{AUT}) are automorphisms is easy to check directly. The argument that proves that this exhausts all automorphisms generated by fields of positive weight is given in \cite{ES98}. It is based on the fact that ${\rm aut} \, \mathcal{Q}$ does not contain fields of weight three and higher. As it turned out \cite{M20}, for an arbitrary nondegenerate $CR$-quadric, this is generally not true. However, there is a simple additional condition satisfied by all $RAQ$-quadrics, which guarantees the validity of this estimate and the statement of item (d) \cite{VB22}.

\vspace{3ex}
Note also that in \cite{WK}, the algebraization procedure was applied to model hyperquadrics. In this case ($k=1$), the transformation $(z \rightarrow z/w, \; w \rightarrow -1/w)$ is a projective automorphism of quadric (\ref{QQ}). This automorphism establishes an isomorphism between $g_j$ and $g_{-j}$ -- the weight components of $\mathrm{aut} \, Q$. This automorphism is preserved under algebraization. Thus, a series of examples of quadratic model surfaces with a nontrivial component $g_{+}$ is constructed.

\vspace{3ex}

{\bf 4. Algebraization of automorphisms}

 \vspace{3ex}

{\bf Statement 5:} Let $Q$ be an arbitrary nondegenerate model surface and $\mathrm{aut}\,Q$ be the Lie algebra of its holomorphic automorphisms. Then there exist two vector-valued polynomials $F$ and $G$ in $(z,w,p)$, where $(z,w)$ are the complex coordinates of the space, and $p=(p_1,\dots,p_d) \in \mathbf{R}^d$ is a set of real parameters, such that $X$ -- an arbitrary element of $\mathrm{aut}\,Q$ -- has the form
\begin{equation}\label{VP}
 X=2 \,\mathrm{Re}\, \left(F(z, w, p) \,\frac{\partial}{\partial z} +G(z,w,p)\, \frac{\partial}{\partial w}\right),
\end{equation}
moreover, $F$ and $G$ depend linearly on $p$, and $d$ is the dimension of $\mathrm{aut}\,Q$.\\
{\it Proof:} $Q$ is nondegenerate, therefore $\mathrm{aut}\,Q$ is finite-dimensional. Let us introduce the grading induced by the Bloom-Graham type. Then each weight component of any field from $\mathrm{aut}\,Q$ is a field from $\mathrm{aut}\,Q$, so $\mathrm{aut}\,Q$ is finitely graded and has the form
 \begin{equation}\label{11}
 2 \,\mathrm{Re}\, \left(\sum_{j=-\nu}^{\mu}\; (F_j(z, w, q) \,\frac{\partial}{\partial z} +G_j(z,w,q)\, \frac{\partial}{\partial w})\right),
 \end{equation}
where $(F_j,G_j)$ are general vector polynomials of the corresponding weights with undetermined real coefficients $q$ in the component of weight $j$. The set of all such coefficients is a finite-dimensional real linear space. The condition for the field (\ref{11}) to be tangent to the surface $Q$ is a system of linear equations with $q$ as unknowns. Its solution is a linear subspace with a system of independent coordinates $p$. The statement is proved.

\vspace{3ex}

Note that if
\begin{eqnarray*}
F(z,w,p)=\sum_{j=1}^{d} \, \varphi_j(z,w)\,p_j, \quad G(z,w,p)=\sum_{j=1}^{d} \, \psi_j(z,w)\,p_j,  \\
X_j=2 \,\mathrm{Re}\, \left(\sum_{j=-\nu}^{\mu}\; (\varphi_j(z, w) \,\frac{\partial}{\partial z} +\psi_j(z,w)\, \frac{\partial}{\partial w})\right),
\end{eqnarray*}
then $(X_1,\dots,X_d)$ is a basis of $\mathrm{aut}\,Q$. The description of $\mathrm{aut}\,Q$ in the form (\ref{VP}) we will call an {\it explicit description}.

We also note that the tangency condition can be divided into a system of independent conditions for each weight component from the $(-\nu)$-th to the $\mu$-th. Therefore, the set of parameters $p$ naturally breaks down into groups of parameters, such that the $j$-th group of parameters (the rest being equal to zero) determines the component $X_j$.

\vspace{3ex}

In order to pass from the Lie algebra $\mathrm{aut}\,Q$ to the local group $\mathrm{Aut}\,Q$, one should solve the following system of first-order ordinary differential equations
\begin{eqnarray}
\nonumber
\frac{d\, \zeta}{d \,t}=F(\zeta,\omega,p),\qquad\qquad\quad \;\qquad \mbox{subject to the condition}\\
\frac{d\, \omega}{d \,t}=G(\zeta,\omega,p),  \qquad\qquad \zeta(0)=z, \quad \omega(0)=w
\end{eqnarray}
with respect to $(\zeta(t),\omega(t))$. If $(\zeta(z,w,p,t),\omega(z,w,p,t))$ is the obtained solution, then
$(\zeta(z,w,p,1),\omega(z,w,p,1))$ is usually denoted as $\mathrm{exp}(X)(z,w,p)$. And we can write that $\mathrm{Aut}\,Q=\mathrm{exp}(\mathrm{aut}\,Q)$.

\vspace{3ex}

Let $Q$ be some nondegenerate model surface of $CR$-type $(n,k)$, let $\mathcal{S}$ be some algebra, $\mathcal{Q}$ be the $\mathcal{S}$-algebraization of $Q$. Let (\ref{VP}) be the explicit description of $\mathrm{aut}\,Q$.
Then the $\mathcal{S}$-{\it algebraization of the Lie algebra} $\mathrm{aut}\,Q$ is the set of $\mathcal{S}$-vector fields of the form
\begin{equation}\label{SVP}
X_{ \mathcal{S}}=2 \,\mathrm{Re}\, \left(F(Z, W, P) \,\frac{\partial}{\partial Z} +G(Z,W,P)\, \frac{\partial}{\partial W}\right),
\end{equation}
where $Z \in \mathcal{S}_c^n, \;\; W \in \mathcal{S}_c^k, \;\; P \in \mathcal{S}^d.$ Notation -- $\mathrm{aut}^s\,Q$.

\vspace{3ex}

Now we can define the $ \mathcal{S}$-{\it algebraization of the local group} $\mathrm{Aut}\,Q$ as
$\mathrm{exp}(\mathrm{aut}_s\,Q)$. Notation -- $\mathrm{Aut}^s\,Q$.

\vspace{3ex}

If we assign some natural weights to the variables $(z_1,\dots,z_n,w_1,\dots,w_k)$, then the spaces of power series and vector fields in $(z,w)$ become graded. Each series and each vector field decomposes into a sum of components of fixed weight.
If during the algebraization $(z,w) \rightarrow (Z,W)$ we transfer the weights to the new algebraic variables, then the same can be said about the spaces of $\mathcal{S}$-power series and $\mathcal{S}$-vector fields, which acquire a grading.

If the variable weights are chosen in accordance with the Bloom-Graham-Stepanova type ({\it induced grading}), then the equations of the model surface $Q$ become weighted homogeneous. This grading extends to vector fields, and the Lie algebra $\mathrm{aut}\,Q$ becomes graded, and in the case of nondegeneracy -- finitely graded. If under algebraization we preserve for the new variables the weights of the corresponding real variables, then the algebraization preserves the graded structure. That is, the algebraization of a weight component is a weight component of the algebraization. The same can be said about the automorphisms, which are the exponential image of a component.
In the graded Lie algebra $\mathrm{aut}\,Q$, one can distinguish the subalgebras $g_{-}$ -- the sum of components of negative weights, $g_0$ -- the component of zero weight, and $g_{+}$ -- the sum of components of positive weights. The exponential maps these subalgebras to local subgroups in $\mathrm{Aut}\,Q$. Notations: $G_{-}, \; G_{0}$, and $G_{+}$.

The algebraization operation commutes with taking the exponential of a weight component. Therefore we obtain:

\vspace{3ex}

{\bf Statement 6:} Let $\mathcal{S}$ be an algebra of dimension $l$. Let $Q$ be a nondegenerate model surface of
$CR$-type $(ln,lk)$ in $\mathcal{S}_c^{(n+k)}$, let $wt$ be the grading corresponding to the type of $Q$ according to Bloom-Graham-Stepanova, $\mathrm{aut}\,Q=\sum_{j=-\nu}^{\mu}\,g_j$ be the decomposition of the Lie algebra of automorphisms corresponding to this grading. Then $\sum_{j=-\nu}^{\mu}\,g^s_j$ is the $wt$-decomposition of $\mathrm{aut}^s\,Q$, with $\mathrm{dim} \,g^s_j=l\,\mathrm{dim} \,g_j$. Respectively, $\mathrm{dim}\,\mathrm{Aut}^s\,Q=l \, \mathrm{dim}\,\mathrm{Aut} \,Q $.

\vspace{3ex}

Thus, under algebraization, the Lie algebra of holomorphic automorphisms only becomes richer. First, algebraization
preserves every nontrivial weight component, mapping it into the corresponding component of the subalgebra of $\mathcal{S}$-holomorphic automorphisms. Second, the dimension of each such component is multiplied by the dimension of $\mathcal{S}$. And third, new holomorphic, but not $\mathcal{S}$-holomorphic automorphisms may appear.

\vspace{1ex}

Let $\mathcal{Q}_{\mathcal{S}}$ be the algebraization of the model surface $Q$.
Let $\mathrm{aut}\,Q$ be the Lie algebra of its holomorphic automorphisms, and $\mathcal{S}\mathrm{aut}\,Q$ be its algebraization.
Is it true that $\mathcal{S}\mathrm{aut}\,Q=\mathrm{aut}\,\mathcal{Q}$? That is, are all fields on the algebraized surface $\mathcal{Q}$ algebraized fields from the original surface $Q$?
The answer to this question is, generally speaking, negative already for the sphere in $\mathbf{C}^2$ (see (\ref{SP})). The Lie algebra of holomorphic fields on the sphere consists of five components with weights from $(-2)$ to $2$. In \cite{ES98} it was shown that the algebraization of the Lie subalgebra consisting of fields of weights 1 and 2 always coincides with the corresponding subalgebra for the algebraization of the sphere. It is easy to show that an analogous statement holds for the subalgebra of fields with negative weights. But for the subalgebra of fields of zero weight, this is not always true.

As an example, one should consider the 2-dimensional algebra generated by the identity and an element $n$ such that $n^2=0$ (dual numbers).
\begin{equation}\label{S2}
S=\{X=\alpha \, 1 + \beta \, n\}, \quad \alpha,\beta \in \mathbf{R}, \quad n^2=0 .
\end{equation}
For the sphere, the Lie subalgebra corresponding to fields of zero weight is 2-dimensional. Therefore, its $\mathcal{S}$-algebraization is 4-dimensional, whereas for
$$\mathcal{Q}_{\mathcal{S}}=\{(z_1, z_2,w_1,w_2) \in \mathbf{C^4}: \mathrm{Im}\, w_1=z_1 \, \bar{z}_1, \;
\mathrm{Im}\, w_2=2 \,\mathrm{Re}\, (z_2 \, \bar{z}_1)\}$$
this subalgebra is 5-dimensional.

Let us find this additional holomorphic but not $\mathcal{S}$-holomorphic field of zero weight. Such a field
has the form
\begin{eqnarray*}
 2 \,\mathrm{Re} \, \left(\lambda(Z)\,\frac{\partial}{\partial Z} +\rho(W)\, \frac{\partial}{\partial W}\right), \; \mbox{where}\;
 Z=z_1 \,1+z_2 \, n, \; W=w_1\,1+w_2\, n, \\
 \lambda(Z)=z_1 \,A_1+z_2 \,A_2, \;\rho(W)=w_1\,B_1+w_2\,B_2, \; A_1, \, A_2 \in \mathcal{S}_c, \; B_1, \, B_2 \in \mathcal{S}.
 \end{eqnarray*}
\begin{eqnarray} \label{TT}
 \mbox{  moreover} \qquad 2 \,\mathrm{Re} \,\left(\lambda(Z)\,\bar{Z}\right)=\rho(Z\,\bar{Z}). \qquad\qquad\qquad\qquad\qquad\qquad
 \end{eqnarray}
For an $\mathcal{S}$-holomorphic field $\lambda(Z)=Z \,A$. Therefore, to find the additional field,
we can pass from an arbitrary field with coefficient $\lambda(Z)$ to $\lambda(Z)-Z\, A_1$. Then $\lambda(Z)=z_2 \,(\alpha+\beta\,n)$. Substituting this into relation (\ref{TT}), we obtain that $\alpha=0, \; \lambda(Z)=z_2 \,\beta \,n, \; \rho(W)=w_2\,\beta\,n, \; \beta \in \mathbf{R}$.

Note that for the remaining 2-dimensional algebras (there are two more of them), this phenomenon is not observed (i.e., all fields are $\mathcal{S}$-holomorphic).

Here is a simple and obvious criterion that all automorphisms are exhausted by $\mathcal{S}$-automorphisms.

\vspace{1ex}

{\bf Statement 7:} All holomorphic automorphisms of the algebraized model surface $\mathcal{Q}_{\mathcal{S}}$ are exhausted by $\mathcal{S}$-automorphisms if and only if the dimension of each weight component of the Lie algebra of automorphisms of the algebraized surface is equal to the dimension of the corresponding component for the original surface $Q$, multiplied by the dimension of the algebra $\mathcal{S}$.

\vspace{3ex}

Let us formulate several questions.

\vspace{1ex}

{\bf Question 8:} Is it true that the algebraization of a nondegenerate germ is a nondegenerate germ?
This question contains two subquestions.\\
(8.1) Is it true that the algebraization of a germ of finite type has finite type?\\
(8.2) Is it true that the algebraization of a holomorphically nondegenerate real-analytic $CR$-manifold is a holomorphically nondegenerate manifold?
\vspace{1ex}

{\bf Question 9:} For which algebras $\mathcal{S}$, model surfaces $Q$, and for which weight components is it true that {\it all} holomorphic fields on $\mathcal{Q}_{\mathcal{S}}$ of a given weight are $\mathcal{S}$-holomorphic?

\vspace{1ex}

{\bf Question 10:} Can $\mathrm{aut}\,\mathcal{Q}$ contain a nonzero component of some weight, provided that there is no component of such a weight in $\mathrm{aut}\,Q$?

\vspace{1ex}

Note that a negative answer to Question 10 ("it cannot") implies a positive answer to Question 8 ("yes, it is nondegenerate").

\vspace{5ex}

{\bf  5. Model surface of degree four}

\vspace{3ex}

In the paper \cite{VB97}, a real hypersurface $Q$ in $\mathbf{C}^2$ with coordinates $(z,w=u+i\,v)$ was considered:
\begin{equation}\label{4}
  Q=\{v=z^2 \, \bar{z}^2\}.
\end{equation}
It was shown there that after introducing the following weights $[z]=[\bar{z}]=1, \; [w]=[\bar{w}]=4$, the algebra of automorphisms ${\rm aut} \, Q$ obtains the following decomposition into weight components $g_{-4}+g_0+g_4$, where
\begin{eqnarray*}
g_{-4}=2 \, {\rm Re}\,( \alpha \,\frac{\partial}{\partial w}), \; \alpha \in \mathbf{R},\\
g_0=2 \, {\rm Re}\,(\beta\, z \,\frac{\partial}{\partial z} +4 \, {\rm Re} ( \beta )\,w \, \frac{\partial}{\partial w}), \;  \beta \in \mathbf{C},\\
g_4=2 \, {\rm Re}\,(\gamma \, z\, w \,\frac{\partial}{\partial z} +2 \,\gamma \,w^2 \, \frac{\partial}{\partial w}), \; \gamma \in \mathbf{R}.
\end{eqnarray*}
The fields of each component generate automorphisms from ${\rm Aut} \, Q$.
\begin{eqnarray*}
g_{-4} \rightarrow G_{-4}= \{ z \rightarrow z,\;\; w \rightarrow w +a\}, \;\; a \in \mathbf{R},\\
g_{0} \rightarrow G_{0}= \{ z \rightarrow b \, z, \;\; w \rightarrow |b|^4 w \}, \;\; b \in \mathbf{C}^{*},\\
g_{4} \rightarrow G_{4}= \{ z \rightarrow \frac{z}{\sqrt{1-c \,w}}, \;\; w \rightarrow  \frac{w}{\sqrt{1-c \,w}}\}, \;\; c \in \mathbf{R}.
\end{eqnarray*}
Consider the algebraization of $Q$.
\begin{equation}\label{Q4}
 \mathcal{Q}_{\mathcal{S}}=\{(Z,W) \in \mathcal{S}_{c} \oplus \mathcal{S}_{c}:  {\rm Im}\,W=Z^2 \,\bar{Z}^2 \}.
\end{equation}
The algebraization of the automorphism group has the form
\begin{eqnarray}  \label{AQ}
 G_{-4}= \{ Z \rightarrow Z,\;\; W \rightarrow W +A\}, \;\; A \in \mathcal{S},\\
\nonumber G_{0}= \{ Z \rightarrow B \, Z, \;\; W \rightarrow B^2 \,\bar{B}^2 \, W \}, \;\; B \in \mathcal{S}_c\\
\nonumber  G_{4}= \{ Z \rightarrow Z \,(1-C \,W)^{(-1/2)}, \;\; W \rightarrow  W \,(1-C \,W)^{(-1/2)}\}, \;\; C \in \mathcal{S}.
\end{eqnarray}
By $(1-C \,W)^{(-1/2)}$ we mean the sum of the corresponding binomial series converging in the norm of $\mathcal{S}_c$, and the invertibility of $B$ is also assumed.

These transformations generate $\mathrm{Aut}_{\mathcal{S}}\,\mathcal{Q}$, a subgroup in $\mathrm{Aut}\,\mathcal{Q}$.
How to describe the entire group $\mathrm{Aut}\,\mathcal{Q}$? For this purpose, consider a general multidimensional analog of (\ref{4}) -- a model surface $Q$ in $\mathbf{C}^{(n+k)}$ of the form
\begin{equation}\label{Q4}
 Q=\{ \mathrm{Im} \, w = \Phi(z,z,\bar{z},\bar{z})\}, \quad z \in \mathbf{C}^n,  \; w \in \mathbf{C}^k, \; n \geq 1, \; k \geq 1
\end{equation}
where $\Phi(z,z,\bar{z},\bar{z})$ is a form of bidegree $(2,2)$ taking values in $\mathbf{R}^k$. Let this surface be nondegenerate, i.e., the coordinates of the form $\Phi$ are linearly independent, and $Q$ is holomorphically nondegenerate.
If
$$ X=\mathrm{Re}\, \left(f(z, w) \,\frac{\partial}{\partial z} +g(z,w)\, \frac{\partial}{\partial w}\right)$$
is a holomorphic vector field in a neighborhood of the origin, then the condition that $X \in {\rm aut} \, Q$ (tangency condition) has the form:
\begin{equation}\label{TC}
 E(z,\bar{z}, u)=2 \, \mathrm{Re}\, (i\,g(z,w)+4\,\Phi(f(z,w),z,\bar{z},\bar{z}))=0, \; w=u+i \, \Phi(z,z,\bar{z},\bar{z}).
\end{equation}
Denote $f(0,u)=a(u), \; g(0,u)=b(u), \; \partial_z (f)(0,u)=c(u)$. Writing the condition $E(z,0,u)=0$ (i.e., setting $\bar{z}=0$ in (\ref{TC})), we obtain $E(z,0,u)=i \, g(z,u)- i \, \bar{b}(u)=0$ and $ g(z,u)= \bar{b}(u)$. Substituting $z=0$, we obtain $\mathrm{Im}\, b(u)=0$ and $ g(z,u)= b(u)$.
Substitute this value for $g$ into (\ref{TC}), differentiate with respect to $\bar{z}$ and substitute $\bar{z}=0$ (i.e., we extract the component linear in $\bar{z}$). We obtain $\bar{\Phi}(\bar{a}(u),\bar{z},z,z)=0$. From which (holomorphic nondegeneracy) we obtain that $a(u)=0$. Now, extracting the component quadratic in $\bar{z}$, we obtain
\begin{eqnarray*}
\partial_u (b)(\Phi(z,z,\bar{z},\bar{z}))=2 \, \Phi(f(z,u),z,\bar{z},\bar{z})+2 \,\bar{\Phi}(\bar{c}(u) \bar{z},\bar{z},z,z).
\end{eqnarray*}
If we write $f(z,u)$ in the form $c(u)\, z+ \varphi(z,u)$, where $\varphi$ are terms that have degree two or higher in $z$, then this relation splits into two (the first is of degree two in $\bar{z}$, the second is higher)
\begin{eqnarray*}\
\partial_u (b)(\Phi(z,z,\bar{z},\bar{z}))=4 \, \mathrm{Re}(\Phi(c(u)\,z ,z,\bar{z},\bar{z})), \quad 2 \,\bar{\Phi}(\bar{\varphi}(\bar{z},u),\bar{z},z,z)=0.
\end{eqnarray*}
If there is a $u=u_0$ such that $ \varphi(z,u_0) \neq 0$, then the field $\varphi(z,u_0) \frac{\partial}{\partial \,z}$ is a nonzero field of type $(1,0)$ tangent to $Q$, which contradicts holomorphic nondegeneracy. Thus, $\varphi=0, \;\; f(z,u)=c(u) z$. The field $X$ has the form
\begin{equation}\label{X}
   X=\mathrm{Re}\, \left(c(w) \,z \, \frac{\partial}{\partial z} +b(w)\, \frac{\partial}{\partial w}\right)
\end{equation}
and the tangency condition takes the form
\begin{equation}\label{EE}
-\mathrm{Im} \, b(u+i\,\Phi(z,z,\bar{z},\bar{z}))+4 \, \mathrm{Re}(\Phi(c(u+i\,\Phi(z,z,\bar{z},\bar{z}))\,z ,z,\bar{z},\bar{z}))=0.
\end{equation}
The series expansion of this expression contains only bidegrees of the form $(2m,2m)$. Denote $\Delta \, \psi(u) =  \partial_u (\psi)(u)(\Phi(z,z,\bar{z},\bar{z}))$. For $m=0,1,2,3$ we obtain
\begin{eqnarray}\label{0123}
\nonumber   \mathrm{Im} \, b(u)=0,  \quad   \Delta(b(u))=4 \, \mathrm{Re}(\Phi(c(u)\,z ,z,\bar{z},\bar{z})),\\
 \Delta(\mathrm{Im} \, \Phi (c(u)\,z ,z,\bar{z},\bar{z}))=0, \quad \Delta^3(b(u))=0.
\end{eqnarray}
The conditions for $m>3$ follow from (\ref{0123}). This system represents a criterion for the field (\ref{X}) to belong to ${\rm aut} \, Q$.

It follows from the theorem of the paper \cite{Z} that, under the condition of nondegeneracy, $\mathrm{deg} \, b \leq 2\,(k+1)$.

Let us summarize our reasoning.

\vspace{3ex}

{\bf Statement 11:} Let $Q$ be a nondegenerate model surface of the form (\ref{Q4}) and $X \in \mathrm{aut}\,Q$, then
(a) $X$ has the form $2\,\mathrm{Re}\, \left(c(w) \,z \, \frac{\partial}{\partial z} +b(w)\, \frac{\partial}{\partial w}\right)$, and the relations (\ref{0123}) are satisfied, (b) $\mathrm{aut}\,Q$ is finite-dimensional, and the degrees of the coefficients do not exceed $2\,(k+1)$.

\vspace{3ex}

Consider the mapping $\mathbf{\Phi}: (z,\zeta) \rightarrow \Phi(z,z,\zeta,\zeta)$. This is a mapping from $\mathbf{C}^{2n}$ to $\mathbf{C}^k$.
Let us formulate a certain condition on the form $\Phi$ ((fd) -- full-dimension).
$$ (fd) \quad \mbox{The image      } \mathbf{\Phi}(\mathbf{C}^{2n})  \mbox{    has interior points.} \qquad \qquad\qquad $$
Note that this condition cannot be satisfied for $k > 2n$.

\vspace{3ex}

The following lemma is quite obvious.

\vspace{1ex}

{\bf Lemma 12:} Under the $(fd)$-condition, $\Delta^m \psi(u)=0$ implies \\ $d^m \psi(u)=0$.

From this lemma we obtain:

\vspace{3ex}

{\bf Statement 13:} If a surface $Q$ of the form (\ref{Q4}) is nondegenerate and $\Phi$ satisfies the $(fd)$-condition, then fields from ${\rm aut} \, Q$ have the form
\begin{eqnarray}\label{Phi1}
 X=2\, \mathrm{Re}\, \left((\lambda(z)+B(w,z))  \, \frac{\partial}{\partial z} +(q + \rho(w)+r(w,w))\, \frac{\partial}{\partial w}\right),\qquad\qquad
\end{eqnarray}
\begin{eqnarray}\label{Phi2}
\nonumber \mbox{moreover  } q, \; \rho \, u, \; r(u,u) \; \mbox{are real and the following relations hold}\qquad,\\
\nonumber 4 \, \mathrm{Re} ( \Phi (\lambda \,z ,z,\bar{z},\bar{z}))=\rho \, \Phi (z ,z,\bar{z},\bar{z}),\qquad\qquad\qquad\qquad \\
\nonumber r(\Phi (z ,z,\bar{z},\bar{z}),u)=2 \, \mathrm{Re} (\Phi (B(u)\,z ,z,\bar{z},\bar{z})), \quad\quad\qquad\qquad\\
\mathrm{Im}\,(\Phi (B(\Phi (z ,z,\bar{z},\bar{z})))\,z ,z,\bar{z},\bar{z}))=0 \qquad\qquad\qquad\qquad
\end{eqnarray}
 {\it Proof:} When this condition is satisfied, it follows from $\Delta^3(b(u))=0$ that $d_u^3(b(u))=0$, i.e., $b(u)=q + \rho (u)+r(u,u)$. Now note that from the second and third relations of (\ref{0123}) it follows that $\Delta^2 \, \Phi (c(u)\,z ,z,\bar{z},\bar{z})=0$, and then $\Phi (d^2(c(u))\,z ,z,\bar{z},\bar{z})=0$. If at the same time $d^2(c(u)) \neq 0$, then the condition of holomorphic nondegeneracy is violated. Thus, $c(u) \,z=(\lambda  + B(u)) \,z$. This determines the form of $X$. Substituting the obtained $b(u)$ and $c(u)$ into (\ref{0123}), we obtain the relations for the parameters. The statement is proved.

\vspace{3ex}

The following statement gives a description of the algebra of all automorphisms of the algebraization of surface (\ref{4}).

\vspace{3ex}

{\bf Statement 14:} If $\mathcal{Q}$ is the algebraization of the surface (\ref{4}), then\\
(a) $\mathcal{Q}$ is nondegenerate,\\
(b) it satisfies the $(fd)$-condition, \\
(c) any field $X \in \mathrm{aut}\,\mathcal{Q}$ has the form (\ref{Phi1}), moreover
\begin{eqnarray}\label{44}
\nonumber 4 \, \mathrm{Re} ( \lambda(Z) Z \bar{Z}^2)=\rho (Z^2 \bar{Z}^2),\qquad\qquad\qquad\qquad  \\
\nonumber 2 \, \mathrm{Re} (B(U,Z) Z \bar{Z}^2)=r(U, Z^2 \bar{Z}^2), \quad\quad\qquad\qquad \\
\mathrm{Im}\,(B(Z^2 \bar{Z}^2, Z) Z \bar{Z}^2)=0 \qquad\qquad\qquad\qquad
\end{eqnarray}
 {\it Proof:}
 Let us prove holomorphic nondegeneracy. Let
 $X=f(Z,W))  \, \frac{\partial}{\partial Z} +g(Z,W)) \, \frac{\partial}{\partial W}$ be a tangent holomorphic $(1,0)$-field.
 Then on $\mathcal{Q}$ we have $g(Z,W)=4 \, i \, Z \,\bar{Z}^2$. Apply the operator
 $ \frac{\partial}{\partial \bar{Z}} - 4\,i \,Z^2\,\bar{Z}\, \frac{\partial}{\partial \bar{W}}$ to this relation, we obtain $f=0$ and $g=0$.
 Holomorphic nondegeneracy is proved. Let us prove the finiteness of type. Let $\alpha(U)$ be a linear functional such that
 $\alpha(Z^2 \, \bar{Z}^2)$ is identically equal to zero. Set $Z=\sqrt{U} \, \bar{Z}=1$, we obtain $\alpha(U)=0$. Item (a) is proved. Substitute the values $A=\sqrt{Z}, \;B=1$ into $\mathbf{\Phi}(A,B)=A^2 \, B^2$, we obtain $\mathbf{\Phi}(A,B)=Z$. This implies the $(fd)$-condition. This is item (b). Now item (c) immediately follows from Statement 13. Statement 14 is proved.

\vspace{3ex}

Let us calculate $\mathrm{aut} \, \mathcal{Q}$ for $Q=\{\mathrm{Im}\,w=z^2\,\bar{z}^2\}$ algebraized using
$\mathcal{S}$ of the form (\ref{S2}). If $(Z=z_1+z_2\,n, \; W=w_1+w_2\,n)=(z_1,z_2,w_1,w_2)$ are coordinates in $\mathcal{S}_c \times \mathcal{S}_c \simeq \mathbf{C}^4$, then the equations for $\mathcal{Q}$ have the form
$$\mathrm{Im}\,w_1=|z_1|^4, \quad  \mathrm{Im}\,w_2=4\,\mathrm{Re}\,(z_1^2 \,\bar{z}_1\,\bar{z}_2).$$
In accordance with Statements 13 and 14, the Lie algebra consists of components of weights (-4), 0, and 4, moreover
$g_{-4}= 2\, \mathrm{Re}\, \left(q \frac{\partial}{\partial W}\right), \; q \in \mathcal{S}$. That is, this component coincides with the algebraization of the corresponding component for the original surface (independent of the algebra).

The field of weight zero has the form
\begin{eqnarray*}
 2 \,\mathrm{Re} \, \left(\lambda(Z)\,\frac{\partial}{\partial Z} +\rho(W)\, \frac{\partial}{\partial W}\right), \; \mbox{where}\;
 Z=z_1 \,1+z_2 \, n, \; W=w_1\,1+w_2\, n, \\
 \lambda(Z)=z_1 \,L_1+z_2 \,L_2, \;\rho(W)=w_1\,r_1+w_2\,r_2, \; L_1, \, L_2 \in \mathcal{S}_c, \; r_1, \, r_2 \in \mathcal{S}\\
 \mbox{  moreover} \qquad 2 \,\mathrm{Re} \,\left(\lambda(Z)\,Z\, \bar{Z}^2\right)=\rho(Z^2\,\bar{Z}^2). \qquad\qquad\qquad\qquad\qquad\qquad
 \end{eqnarray*}
Using the fact that we have the $\mathcal{S}$-holomorphic field $\lambda(Z)=L \, Z, \; \rho(W)=2\,\mathrm{Re}(L) \,W$, we set $L_1=0$. Then substituting into the relation we obtain $r_1=0, \; L_2=2 \, \alpha \, n, \;r_2= \alpha \, n, \; \alpha \in \mathbf{R}$. That is, the subspace of non-$\mathcal{S}$-holomorphic fields in $g_0$ is 1-dimensional and is generated by the field
$$2\, \mathrm{Re}\,
(z_2 \, \frac{\partial}{\partial z_2} +2 \,w_2 \, \frac{\partial}{\partial w_2})
$$
The $\mathcal{S}$-holomorphic part of $g_0$ is 4-dimensional, and the entire $g_0$ is 5-dimensional.

\vspace{3ex}

Proceeding to the calculation of $g_4$, we set
\begin{eqnarray*}
B(U,Z)=u_1 \, (B_{11}\,z_1+B_{12}\,z_2)+u_2 \, (B_{21}\,z_1+B_{22}\,z_2), \\
r(U,V)=u_1 \, (r_{1}\,v_1+r_{2}\,v_2)+u_2 \, (r_{2}\,v_1+r_{3}\,v_2),\\
B_{11},\, B_{12},\,B_{21},\,B_{22} \in \mathcal{S}_c, \; r_{1},\,r_{2},\,r_{3} \in \mathcal{S},.
\end{eqnarray*}
Taking into account that the $\mathcal{S}$-holomorphic part of $g_4$ contains the field $(R\,Z\,W, \;2\,R\,W^2)$, where $R \in \mathcal{S}$, we can assume that $r_1=0$. Then substituting these expressions into the second and third relations of (\ref{44}) and separating the coefficients of different monomials, we obtain that the only solution of this form is $B(U,Z)=0, \;  r(U,V)=0$.
That is, there are no fields in $g_4$ that are not $\mathcal{S}$-holomorphic. Thus, we can summarize:

\vspace{3ex}

{\bf Statement 15:} Let $\mathcal{S}$ be a 2-dimensional algebra consisting of elements of the form $X=x_1\,1+x_2\,n$, where $n^2=0$.
Let $\mathcal{Q} =\{\mathrm{Im\,W}=Z^2\,\bar{Z}^2\}, \; (Z,W) \in \mathcal{S}^2_c$. Then the complement of $\mathrm{aut}_s\,\mathcal{Q}$ to $\mathrm{aut} \,\mathcal{Q}$ is a 1-dimensional space generated by the field
 $2\, \mathrm{Re}\,(z_2 \, \frac{\partial}{\partial z_2} +2 \,w_2 \, \frac{\partial}{\partial w_2})$. It corresponds to a 1-parameter group of automorphisms of the form
 $$(z_1 \rightarrow z_1, \; z_2 \rightarrow e^t \,z_2, \; w_1 \rightarrow w_1, \; w_2 \rightarrow e^{2t} \,w_2),$$
whereas the $\mathcal{S}$-holomorphic part of the local group has the form (\ref{AQ}).

\vspace{3ex}

{\bf 6. Concluding remarks}

\vspace{3ex}

{\bf(A)} Algebraization in the context of $CR$-geometry is, on the one hand, a way of multiplying interesting examples, and on the other, a procedure for identifying different algebraizations of the same original object. In this regard, whenever a highly symmetric object arises, one should try to understand it as the result of the algebraization of an already known, simpler object.

\vspace{1ex}

There are many interesting examples for replication. Among them:
light cones in $\mathbf{C}^{n+1}$ -- $\{(\mathrm{Im}\,z_{n+1})^2=\sum_{j=1}^{n} \, (\mathrm{Im}\,z_{j})^2\}$,
lists of 3-, 4-, 5-dimensional holomorphically homogeneous manifolds \cite{BK}, \cite{LB}, series of Labovsky hypersurfaces \cite{LA}, Zelenko and Sykes \cite{ZS}, examples of exceptional quadrics by Meylan and Gregorovi$\check{c}$ \cite{GM}, 3-nondegenerate Santi hypersurfaces \cite{SN}, extremely symmetric analytic hypersurfaces of infinite type in $\mathbf{C}^2$ with 5-dimensional stabilizers \cite{VB12}
$$v=u\, \mathrm{tg}\left(\frac{1}{2m}\,\mathrm{arcsin}(2m\,z\,\bar{z})  \right),   \quad m=1,2,\dots$$

\vspace{3ex}

{\bf(B)} Some possibilities of algebraization in $CR$-geometry were demonstrated above. Moreover, this was a demonstration within the framework of the {\it analytic} point of view (coordinates, equations, ...), i.e., $CR$-geometry $\grave{a}$ la Poincar$\acute{e}$. An alternative, {\it geometric} view on $CR$-geometry is $CR$-geometry $\grave{a}$ la Cartan (see \cite{T}). In this case, the base object is a fundamental graded Lie algebra.
This approach also allows algebraization: we algebraize the fundamental Lie algebra, obtain an $\mathcal{S}$-fundamental graded Lie algebra, and then launch the Tanaka prolongation process.

{\bf(C)} The possibilities of algebraization undoubtedly go beyond the scope of $CR$-geometry. This procedure is appropriate in the context of broadly understood {\it analytic geometry}. For example, there are known lists \cite{MHZ} of homogeneous affine manifolds of small dimensions. They can all be algebraized, and homogeneity will be preserved.
The same applies to projectively homogeneous manifolds.
All constructions of linear algebra and projective geometry admit a natural algebraization.
Algebraization is, of course, applicable to algebraic manifolds. Undoubtedly, algebraists are aware of this possibility. Algebraic groups can be subjected to algebraization. If we have a linear representation of some Lie algebra, it is appropriate to consider the algebraization of this Lie algebra by algebraizing the representation space.
The entire mathematical analysis that does not leave the boundaries of the {\it analytic universe} can be subjected to algebraization. All constructions of complex analysis (power series, Cauchy's integral formula) are naturally algebraized. One can speak about $\mathcal{S}$-entire and $\mathcal{S}$-meromorphic functions. In a sense, we have as many copies of analysis as we have algebras.
The same applies to the analytic theory of differential equations.
If we have algebraized an ordinary analytic differential equation, we obtain its $\mathcal{S}$-analog, the solutions of which will be $\mathcal{S}$-holomorphic functions. From a point of view external to algebraization, such a system is a system of partial differential equations.
Under the algebraization of systems of linear partial differential equations with constant coefficients, a subspace of $\mathcal{S}$-holomorphic solutions appears in the solution space. One can propose an $\mathcal{S}$-version of the Cauchy-Kovalevskaya theorem.
One can algebraize the theory of elliptic functions, both in the Weierstrass version ($\wp$-function) and the Poincar$\acute{e}$ version ($\Theta$-series). This same applies to all special functions, including $\Gamma(z)$ and $\zeta(s)$. The Gauss hypergeometric function $F(a,b;c;z)$ has its own algebraic version. As well as hypergeometric functions of several variables.
It is not difficult to give a definition of an $\mathcal{S}$-analytic and $\mathcal{S}$-complex manifold and to develop a corresponding version of analysis ($\mathcal{S}$-meromorphic functions on an $\mathcal{S}$-complex manifold).
The procedure for algebraizing toric and spherical manifolds also seems interesting and meaningful.

\vspace{3ex}

{\bf(D)} Undoubtedly, one can weaken the requirements for the algebra used to perform the algebraization.
There are quaternions, octonions, Clifford algebras... There are infinite-dimensional algebras.

\vspace{3ex}

{\bf(E)} All this can be understood as the outlines of a certain total program -- a program of $\mathcal{S}$-algebraization of everything in the world. This is not the first general mathematical program and, of course, the expediency of applying the algebraization procedure in a given area is a priori not obvious, but not excluded.

\end{document}